\documentclass[11pt,reqno]{amsart}

\usepackage{color}

\newtheorem{theorem}{Theorem}
\newtheorem*{theoremNN}{Theorem}
\newtheorem{corrolaire}{Corollaire}

\newtheorem{proposition}[theorem]{Proposition}
\newtheorem{corollary}[theorem]{Corollary}
\newtheorem{example}{Example}
\newtheorem{remark}{Remark}
\newtheorem{lemma}[theorem]{Lemma}
\newfont{\bb}{msbm10 at 12pt}

\def\pf{\noindent{\textit {Proof.} }}

\def\Ric{\hbox{Ric}}

\newcommand{\bal}{\begin{align}}      \newcommand{\eal}{\end{align}}
\newcommand{\ba}{\begin{array}}      \newcommand{\ea}{\end{array}}
\newcommand{\bc}{\begin{center}}     \newcommand{\ec}{\end{center}}
\newcommand{\be}{\begin{enumerate}}  \newcommand{\ee}{\end{enumerate}}
\newcommand{\beq}{\begin{eqnarray}}  \newcommand{\eeq}{\end{eqnarray}}
\newcommand{\beQ}{\begin{eqnarray*}} \newcommand{\eeQ}{\end{eqnarray*}}
\newcommand{\bi}{\begin{itemize}}    \newcommand{\ei}{\end{itemize}}
\newcommand{\bt}{\begin{tabular}}    \newcommand{\et}{\end{tabular}}
\newcommand{\bdm}{\begin{displaymath}} \newcommand{\edm}{\end{displaymath}}

\newcommand{\trg}{{\rm Tr}_\gamma}
\newcommand{\trgr}{{\rm Tr}_{g_r}}

\def\qed{\hfill{q.e.d.}\smallskip\smallskip}

\begin{document}

\title[Cheeger constant and conformal infinity]{The Cheeger constant of an asymptotically locally hyperbolic manifold and the Yamabe type of its conformal infinity}  
 
\author{Oussama Hijazi}
\address[Oussama Hijazi]{Institut {\'E}lie Cartan,
Universit{\'e} de Lorraine, Nancy,
B.P. 70239,
54506 Vand\oe uvre-L{\`e}s-Nancy Cedex, France.}
\email{Oussama.Hijazi@univ-lorraine.fr}

\author{Sebasti{\'a}n Montiel}
\address[Sebasti{\'a}n  Montiel]{Departamento de Geometr{\'\i}a y Topolog{\'\i}a,
Universidad de Granada,
18071 Granada,  Spain.}
\email{smontiel@ugr.es}

\author{Simon Raulot}
\address[Simon Raulot]{Laboratoire de Math\'ematiques R. Salem
UMR $6085$ CNRS-Universit\'e de Rouen
Avenue de l'Universit\'e, BP.$12$
Technop\^ole du Madrillet
$76801$ Saint-\'Etienne-du-Rouvray, France.}
\email{simon.raulot@univ-rouen.fr}

\begin{abstract}
Let $(M,g)$ be an $(n+1)$-dimensional asymptotically locally hyperbolic (ALH) manifold with a conformal compactification whose conformal infinity is $(\partial M,[\gamma])$. We will first observe that ${\mathcal Ch}(M,g)\le n$, where  ${\mathcal Ch}(M,g)$ is the Cheeger constant of  $M$.  We then  prove that, if the Ricci curvature of $M$ is bounded from below by $-n$  and its scalar curvature approaches $-n(n+1)$ fast enough at infinity,   then ${\mathcal Ch}(M,g)= n$  if and only ${\mathcal Y}(\partial M,[\gamma])\ge 0$, where ${\mathcal Y}(\partial M,[\gamma])$ denotes the Yamabe invariant of the conformal infinity. This gives 
 an answer to a question raised  by J. Lee \cite{L}.
\end{abstract}



\date {\today}        

\maketitle \pagenumbering{arabic}


\section{Introduction}


We will study asymptotically locally hyperbolic (ALH) manifolds in the more general setting of conformally compact manifolds. 
During the last decades, due to the important 
role that they play in the so-called anti-de Sitter/conformal field theory (AdS/CFT) correspondence (see
\cite{Bi}, for instance), this class of Riemannian manifolds has attracted a great deal of interest 
in both physical and mathematical realms. The mathematical 
aspects relative to the existence and the behavior near the infinity of conformally compact  manifolds satisfying the Einstein condition were first
studied in the seminal paper  by C. Fefferman and C. Graham \cite{FG}. This particular class of conformally compact manifolds are the usually called Poincar{\'e}-Einstein (PE) spaces (they necessarily have negative constant scalar curvature). 

As in many papers about asymptotically hyperbolic manifolds, here we will drop the Einstein condition and call ALH manifold any conformally compact Riemannian manifold whose scalar curvature is asymptotically constant (also necessarily negative). This implies that the same must occur for its sectional curvatures. In some sense,  ALH manifolds are just those  looking like PE spaces at infinity. On the other hand, a Riemannian manifold $(M,g)$ is called conformally compact if it is a connected complete manifold whose metric extends conformally
 to a compact manifold with (non-necessarily connected) boundary whose interior is the original manifold.
So, by means of this extended conformal metric, the corresponding original metric determines a conformal structure on the boundary $(\partial M,[\gamma])$, which is usually called the conformal infinity or the boundary at infinity. 

In this setting, a natural question is to look for relations between Riemannian invariants of an $(n+1)$-dimensional ALH manifold $M$, $n\ge 2$, and  conformal invariants  of the $n$-dimensional conformal boundary $(\partial M,[\gamma])$. A beautiful result in this direction was obtained by J. Lee in \cite{L}. In fact, he proved that if ${\mathcal Y}(\partial M,[\gamma])\ge 0$, then $\lambda_{1,2}(M)=\frac{n^2}{4}$, where ${\mathcal Y}(\partial M,[\gamma])$ denotes the Yamabe invariant of the compact conformal manifold $(\partial M,[\gamma])$, which is the infimum of the total scalar curvature functional over unit-volume metrics in the conformal class $[\gamma]$, and $\lambda_{1,2}(M)$ is the infimum of the $L^2$ spectrum of the Laplacian of $M$. In this way, he thoroughly extended a result which was known to occur when $M={\mathbb H}^{n+1}/\Gamma$ is a geometrically finite and cusp-free quotient of the hyperbolic space by a Kleinian group, a consequence from 
previous works by D. Sullivan and by R. Schoen and S.-.T. Yau (see \cite{Su,SY}). J. Lee pointed
out 
that his theorem is not sharp, in the sense that there are ALH manifolds $M$ with ${\mathcal Y}(\partial M,[\gamma])<0$  but still $\lambda_0(M)=\frac{n^2}{4}$, and raised the question of finding a necessary and sufficient condition on the geometry of $M$ for ${\mathcal Y}(\partial M,[\gamma])\ge 0$. In this direction, C. Guillarmou and J. Qing proved in \cite{GQ} that, when $M$ is a PE space with $n>2$, ${\mathcal Y}(\partial M,[\gamma])> 0$ if and only if  the largest real scattering pole of $M$ is less than $\frac{n}{2}-1$. 

In this work, we answer the aforementioned question (see Theorem \ref{maintheorem}) by relating the Yamabe type of the boundary at infinity with the value of the Cheeger constant ${\mathcal Ch}(M,g)$ of the ALH manifold $(M,g)$ (see Section \ref{SectionUpper} for a precise definition), namely 
\begin{theoremNN}
Let $(M,g)$ be an $(n+1)$-dimensional conformally compact Riemannian manifold of order $C^{m,\alpha}$ with $m\geq 3$, $0<\alpha<1$ and whose Ricci tensor and scalar curvature satisfy
$$
\hbox{\rm Ric}_g+ng\ge 0,\qquad R_g+n(n+1)=o\big(r^2\big),
$$
where $r$ is any defining function on $M$, then 
$$
{\mathcal Ch}(M,g)= n\quad\Longleftrightarrow\quad{\mathcal Y}(\partial M,[\gamma])\ge 0.
$$
\end{theoremNN}
Since it is not difficult to observe that ${\mathcal Ch}(M,g)\le n$ for each $(n+1)$-dimensional ALH manifold (see Corollary \ref{upper-Cheeger}), an equivalent statement of our characterization is to say that
${\mathcal Y}(\partial M,[\gamma])\ge 0$ if and only if  the following linear isoperimetric inequality
\begin{equation}\label{LII}
A(\partial\Omega)\ge nV(\Omega)
\end{equation}
holds for all compact domains $\Omega\subset M$ (see Corollary \ref{yamabe-yau}). This isoperimetric inequality was well-known to be valid for hyperbolic spaces and S.-T. Yau proved that it is also true on complete simply connected Riemannian manifolds with sectional curvatures bounded from above by $-1$  (see \cite{Y} and \cite[Theorem 34.2.6]{BZ}). Another direct consequence of our main result is the generalization of the result of Lee on the bottom of the $L^2$ spectrum of the Laplacian of $M$ to the principal eigenvalue of its $p$-Laplacian (see Theorem \ref{lambda}).

It is worth mentioning (and maybe useful to the reader) that, in a different context of hyperbolicity, J. Cao \cite{C} also explored 
 the relation between the geometric properties of a Gromov-hyperbolic space and some diverse features of its boundary at infinity.


\section{Conformally compact Riemannian manifolds}


Let ${\overline M}$ be a (connected) compact $(n+1)$-manifold with  (non-ne\-ce\-ssa\-ri\-ly connected) boundary and $n\ge 2$. The interior of ${\overline M}$ will be denoted by $M$ and its boundary by $\partial M$. If $g$ is a smooth Riemannian metric on $M$, the open Riemannian manifold $(M,g)$ is said to be {\em conformally compact} of order $C^{m,\alpha}$ if, for some (and hence for all) smooth defining function $\rho$, the smooth conformal metric $\rho^2 g$ on $M$ extends to a $C^{m,\alpha}$ Riemannian metric $\overline{g}$ on ${\overline M}$. Here $C^{m,\alpha}$ denotes the classical H\"older space for $m\in\mathbb{N}$ and $\alpha\in[0,1]$. Recall that a $C^1$ map $\rho:{\overline M}\rightarrow\mathbb{R}$ is a {\em defining function} of the boundary if it is a non-negative function such that $\rho^{-1}(\{0\})=\partial M$ and $d\rho\ne 0$ everywhere on $\partial M$. It is obvious that there are many different defining functions for $\partial M$, but all the corresponding extended metrics $\overline{g}=\rho^2 g$ 
will have conformally equivalent restrictions to the boundary $\partial M$. Then, if $\gamma=\overline{g}_{|\partial M}$, the conformal manifold $(\partial M,[\gamma])$  
is well defined and depends only on $(M,g)$. This pair is called the {\em conformal infinity} of $(M,g)$. 
 
The simplest example of a conformally compact Riemannian manifold is the hyperbolic space ${\mathbb H}^{n+1}$ which can be realized as the Riemannian manifold $\big(B^{n+1},\frac{4|dx|^2}{(1-|x|^2)^2}\big)$. Here $B^{n+1}=\{x\in {\mathbb R}^{n+1}\,/\,|x|<1\}$ is the $(n+1)$-dimensional Euclidean unit open ball and $|dx|^2$ is the flat Euclidean metric. In this situation, the map 
$x\in B^{n+1}\mapsto(1-|x|^2)/2$ is a defining function for the boundary $\partial B^{n+1}={\mathbb S}^n$ and the corresponding conformal infinity is then easily seen to be $({\mathbb S}^n,[g_0])$, where $[g_0]$ denotes the conformal class of the round metric $g_0$ of constant sectional curvature $1$ on ${\mathbb S}^n$.    
 
Assume now that the conformally compact Riemannian manifold $(M,g)$ is of order at least $C^{2}$. Then using the relation of the Riemannian curvature tensors under conformal changes of metrics (see, for instance, \cite[p.\! 59]{Be}), it can be easily seen that 
 all its sectional curvatures $K_g$ uniformly approach 
$-|{\overline \nabla}\rho|_{\overline{g}}^2$ as $\rho\rightarrow 0$. Here ${\overline\nabla}$ is the gradient operator corresponding to the metric ${\overline g}$ and the norm is taken with respect to the same metric as that of the gradient. So, it is clear that the quantity $|{\overline \nabla}\rho|_{\overline{g}}$ restricted to $\partial M$ depends only on the original metric $g$. Thus, {\em conformally compact manifolds of order at least $C^{2}$ are asymptotically negatively curved}. From this observation, we will say that a conformally compact Riemannian manifold $(M,g)$ is {\em asymptotically locally hyperbolic (ALH)} when $|{\overline \nabla}\rho|_{\overline{g}|\partial M}$ is constant, that we normalize to be equal to $1$. It immediately follows that we have $K_g\rightarrow -1$ near infinity and so the Ricci tensor satisfies $\Ric_g\rightarrow-ng$, that is, the manifold seems to be Einstein with Ricci curvature $-n$ when one moves towards infinity. Obviously, the scalar curvature $R_g$ tends to the constant value $-n(n+1)$. 

Conversely, from the transformation rules of the Ricci tensor and the scalar curvature of conformal metrics, it can be seen that any of these asymptotical behaviors for $K_g$, $\Ric_g$ or $R_g$ implies that $|{\overline \nabla}\rho|_{\overline{g}|\partial M}=1$ for any defining function $\rho$. This means that {\em a $C^{2}$ conformally compact Riemannian manifold is ALH if and only if it is asymptotically Einstein, that is, $\hbox{\rm Ric}_g+ng\rightarrow 0$ uniformly}. As we just noticed, it is also equivalent to the fact that {\em the scalar curvature is asymptotically constant, that is, $R_g+n(n+1)\rightarrow 0$ uniformly as $\rho\rightarrow 0$}. In particular, this occurs when the manifold $(M,g)$ is supposed to be Einstein. In this case, $(M,g)$ is often called a Poincar{\'e}-Einstein manifold (in short, a PE manifold) and we have $\Ric_g+ng=0$. The weaker condition on the constant scalar curvature $R_g +n(n+1)=0$ implies that $(M,g)$ is an ALH manifold as well.
 
In the general non-necessarily Einstein ALH case, if we assume $(M,g)$ to be conformally compact of order at least $C^{3,\alpha}$, we can modify any given smooth defining function $\rho$ to get another one $r\in C^{2,\alpha}({\overline M})$ such that the corresponding extended conformal metric ${\overline g}= r^2g$ is of class $C^{2,\alpha}$ and  $|{\overline\nabla r}|_{\overline{g}}\equiv 1$ in a collar neighborhood of the boundary at infinity $\partial M$. More precisely, we have 
\begin{lemma}{\rm ({\cite[Lemma 5.2]{GL},  \cite[Lemma 5.1]{L}}\label{defining-geodesic})}
Let $(M,g)$ be an ALH manifold  of class $C^{m,\alpha}$ with $m\ge 3$ and $0<\alpha<1$. For each choice of a metric $\gamma$ on its conformal infinity $(\partial M,[\gamma])$, there exists a defining function
$r\in C^{m-1,\alpha}({\overline M})$ uniquely determined in a neighborhood of $\partial M$ such that the extended conformal metric ${\overline g}=r^2g$ is of class $C^{m-1,\alpha}$,  $|{\overline\nabla} r|_{\overline{g}}\equiv 1$ in this neighborhood and with
\begin{equation}\label{Taylor}
{\overline g}=dr^2+g_r=dr^2+\gamma+rg^{(1)}+r^2g^{(2)}+O(r^{2+\alpha}),
\end{equation}
where $O(r^{2+\alpha})$ is a symmetric two-tensor on $\partial M$. Moreover $g^{(i)}$ is of class $C^{2-i,\alpha}$ for $i=0,1,2$ and is computable from the iterated Lie derivatives of the extended metric:
\begin{equation}\label{Lie} 
g^{(i)}=\frac1{i!}\left.{\mathcal L}_{\overline\nabla r}^{(i)}{\overline g}\right|_{r=0}.
\end{equation}
We will say that such a function $r$ is {\em the geodesic defining function} associated with the metric $\gamma$.
\end{lemma} 

For these reasons, we will always assume in this paper that the ALH hyperbolic manifolds considered are of class $C^{m,\alpha}$ with $m\geq 3$ in such a way that for any choice of a metric in the conformal infinity we have a compactification of class $C^{2,\alpha}$ for which the compactified metric has an expansion given by (\ref{Taylor}). A particularly interesting class of such manifolds are the PE spaces or, as discussed in the next section, the weakly Poincar\'e-Einstein (WPE) spaces.


\section{Upper bounds for the Cheeger constant of ALH manifolds}\label{SectionUpper}


In this section, we will see that the Yamabe type of the conformal infinity $(\partial M,[\gamma])$ has a direct influence on the isoperimetric behavior of the large regions of $(M,g)$. In fact, we will observe that these properties can be expressed using the Cheeger constant of $M$.

In a first step, we collect some curvature properties for level hypersurfaces (near infinity) of any geodesic defining function of the boundary. More precisely, suppose that $(M,g)$ is an ALH manifold and fix $\gamma\in[\gamma]$ and $r$ the corresponding geodesic defining function. For $r>0$ sufficiently small, the level sets $\Sigma_r=\{r={\rm const}.\}$ are smooth compact embedded hypersurfaces. If $H_r$ denotes the (inner) mean curvature of $\Sigma_r$, it is straightforward to observe from the first equality in (\ref{Taylor}) that 
\begin{eqnarray}\label{MeanCurvatureF1}
H_r =  \frac{1}{2n}r^{2}\trgr\big(-r\partial_r(r^{-2}g_r)\big) = 1-\frac{r}{2n}\trgr(\partial_rg_r).
\end{eqnarray}
Then using the second equality in (\ref{Taylor}), we compute that
\begin{eqnarray}\label{MeanCurvatureF2}
H_r =  1-\frac{r}{2n}\trg(g^{(1)})+\frac{1}{n}\Big(\frac{1}{2}\trg(A_{\overline{g}}^2)-\trg(g^{(2)})\Big)r^2+O(r^{2+\alpha})
\end{eqnarray}
where $A_{\overline{g}}$ is the symmetric endomorphism of the tangent bundle of $\partial M$ with respect to $\gamma$ defined by $A_{\overline{g}}:=\gamma^{-1}g^{(1)}$. We immediately deduce from (\ref{MeanCurvatureF1}) that $H_r$ is $C^{2,\alpha}$ in $r$ and $H_0=1$. 

Assume now for a moment that the manifold is weakly Poincar\'e-Einstein (WPE) in the sense that the coefficients in the asymptotic expansion (\ref{Taylor}) of $g$ are given by 
$g^{(1)}=0$ and
\begin{eqnarray*}
\quad{g^{(2)}}=-P_\gamma=-\frac{1}{n-2}\Big(\hbox{\rm Ric}_\gamma-\frac{R_\gamma}{2(n-1)}\gamma\Big)
\end{eqnarray*} 
where $\hbox{\rm Ric}_\gamma$, $R_\gamma$ and $P_\gamma$ denote respectively the Ricci tensor, the scalar curvature and the Schouten tensor of $\gamma$. It can be shown that these conditions are equivalent to a second order decay assumption of the Ricci tensor of the metric $g$ namely $|\hbox{\rm Ric}_g+ng|_g=o(r^2)$. Then in this situation we observe that formula (\ref{MeanCurvatureF2}) directly implies that $H'_0=0$ and $H''_0=R_\gamma/(n(n-1))$ so that the value on the boundary of the second derivative of the mean curvature of the level sets of the geodesic defining function with respect to $\gamma$ is precisely encoded by the scalar curvature of this metric. The purpose of the next proposition is to show that a similar result holds under weaker curvature assumptions. By keeping the notations introduced in the above discussion, we have 
\begin{proposition}\label{H-limit} 
 Let $(M,g)$ be an $(n+1)$-dimensional ALH manifold of order $C^{m,\alpha}$ with $m\geq 3$ and $0<\alpha<1$. Then $H_r$ extends to a $C^{2,\alpha}$ function at $r=0$ with $H_0=1$. If, in addition, we have $R_g+n(n+1)=o\big(r^2\big)$, then 
\begin{equation}\label{meancur-1}
H'_0=0, \qquad H''_0\le \frac{R_\gamma}
{n(n-1)}.
\end{equation}
If, moreover, we have $\hbox{\em Ric}_g\ge -ng$, then for $\varepsilon>0$ sufficiently small,
\begin{equation}\label{meancur-2}
H_r\ge  1+\frac{R_\gamma}{2n(n-1)}r^2,\qquad
0\le r\le \varepsilon.
\end{equation}     
\end{proposition}

\begin{remark}{\rm 
It is clear that the decay conditions on the scalar and the Ricci curvatures are slightly hardening  the ALH condition. Note however that they are obviously satisfied when $(M,g)$ is a PE space or a WPE space. 
}\end{remark}

\pf 
As before, we will work in a collar neighborhood of $\partial M$ where $|{\overline \nabla}r|_{\overline{g}}^2= 1$ for $r$ the unique geodesic defining function associated with $\gamma$ whose existence is ensured by Lemma \ref{defining-geodesic}. Note that we have already proved that $H_r$ is a $C^{2,\alpha}$ function in $r$ with $H_0=1$. 

First observe that since $g$ and $\overline{g}=r^2g$ are two conformally related metrics, their scalar curvatures satisfy 
\begin{equation}\label{step2-r}
\frac1{r}\big(R_g+n(n+1)\big)=rR_{\overline g}+2n{\overline\Delta}r.
\end{equation}
So if we assume that $R_g+n(n+1)=o(r^2)$, the previous identity implies in particular that ${\overline\Delta}r_{|\partial M}=0$. On the other hand, since $\overline{g}=dr^2+g_r$ we compute that
\begin{eqnarray}\label{LaplaceForm-1}
{\overline\Delta}r_{|\partial M}=\frac{1}{2}\trg(g^{(1)})
\end{eqnarray}
and then $H'_0=-\trg(g^{(1)})/(2n)=0$ where the first equality follows from (\ref{MeanCurvatureF2}). Moreover, since
\begin{eqnarray}\label{VolumeForm}
\frac{\partial_r\sqrt{{\rm det}(g_r)}}{\sqrt{{\rm det}(g_r)}}=\frac{1}{2}\trgr(\partial_rg_r)=r\Big(\trg(g^{(2)})-\frac{1}{2}\trg(A_{\overline{g}}^2)\Big)+O(r^{1+\alpha})
\end{eqnarray}
we also deduce from (\ref{step2-r}) that
\begin{eqnarray}
\frac{R_g+n(n+1)}{r^2}=R_{\overline{g}}+2n\Big(\trg(g^{(2)})-\frac{1}{2}\trg(A_{\overline{g}}^2)\Big)+O(r^{\alpha}).
\end{eqnarray}
Our assumption on the asymptotic behavior of the scalar curvature implies that 
\begin{eqnarray}\label{DerSecMean}
H''_0=-\frac{2}{n}\Big(\trg(g^{(2)})-\frac{1}{2}\trg(A_{\overline{g}}^2)\Big)=\frac{1}{n^2}R_{\overline{g}|\partial M}
\end{eqnarray}
where the first equality follows from (\ref{MeanCurvatureF2}). 

Now we note that the mean curvature $H_r$ is easily computable using the well-known relation between the two mean curvatures of a hypersurface corresponding to two conformal metrics $g=\frac1{r^2}{\overline g}$ on the ambient space (see \cite{E}, for instance):
\begin{eqnarray}\label{step4-r}
H_r=r\big({\overline H}_r-{\overline g}( {\overline\nabla}\log \frac1{r},{\overline N}_r)\big)=1+r{\overline H}_r.
\end{eqnarray} 
Here ${\overline N}_r=\overline{\nabla}r$ (resp. ${\overline H}_r$) denotes the inner unit normal (resp. the mean curvature) of $\Sigma_r$ with respect to the metric ${\overline g}$.  This identity immediately implies that ${\overline H}_0=0$.

Moreover, since $r$ is in fact the $\overline{g}$-distance from the boundary, the mean curvature function $\overline{H}_r$ satisfies the well-known Riccati equation (which can be deduced from \cite[p.\! 44]{P})
\begin{eqnarray}\label{bochner}
n{\overline H}'_r=|{\overline\sigma}_r|_{\overline{g}}^2
+\hbox{\rm Ric}_{\overline g}(\overline{N}_r,\overline{N}_r)
\end{eqnarray}
where $\overline{\sigma}_r$ is the second fundamental form of $\Sigma_r$ with respect to the metric ${\overline g}$. On the other hand, since ${\overline H}_0=0$, the Gau{\ss} formula implies that
$$
\hbox{\rm Ric}_{\overline g}({\overline N}_0,{\overline N}_0) = \frac{1}{2}\Big(R_{{\overline g}|\partial M}-R_\gamma-|{\overline\sigma}_0|_{\overline{g}}^2\Big)
$$
and so (\ref{bochner}) for $r=0$ writes
\begin{eqnarray*}
\overline{H}'_0=\frac{1}{2n}\left({R_{\overline g}}_{|\partial M}-R_\gamma+|{\overline\sigma}_0|_{\overline{g}}^2\right).
\end{eqnarray*}
This with formula (\ref{step4-r}) yields to
\begin{eqnarray*}
H''_0=\frac{1}{n}\left({R_{\overline g}}_{|\partial M}-R_\gamma+|{\overline\sigma}_0|_{\overline{g}}^2\right)
\end{eqnarray*}
which, when combined with (\ref{DerSecMean}), gives
\begin{eqnarray*}
H''_0=\frac{1}{n(n-1)}\left(R_\gamma-|{\overline\sigma}_0|_{\overline{g}}^2\right).
\end{eqnarray*}
The inequality in (\ref{meancur-1}) follows directly. 

We assume now in addition that $\hbox{\rm Ric}_g\ge -ng$. Since $g$ and $\overline{g}$ are conformally related we compute (see \cite[p.\! 59]{Be}) that 
\begin{equation}\label{step1-r}
r(\hbox{\rm Ric}_g+ng)=r\hbox{\rm Ric}_{\overline g}+(n-1) {\overline\nabla}^2r+({\overline\Delta}r) {\overline g}
\end{equation} 
where ${\overline\nabla}^2$ denotes the Hessian of a function with respect to $\overline{g}$. Applying this formula to the vector field $\frac{{\overline\nabla}r}{r}$ and using the fact that $|\overline{\nabla} r|_{\overline{g}}=1$ yield to
\begin{eqnarray}\label{ricn-1}
\hbox{\rm Ric}_{\overline g}({\overline N}_r,{\overline N}_r)-\frac{n\overline{H}_r}{r}\geq 0
\end{eqnarray}
since ${\overline H}_r=-\frac1{n}{\overline\Delta}r$. Taking the limit as $r\rightarrow 0$ implies that 
\begin{equation}\label{ricn-2}
\hbox{\rm Ric}_{\overline g}({\overline N}_0,{\overline N}_0)\geq n\overline{H}'_0.
\end{equation}
However from (\ref{bochner}) with $r=0$, we observe that this inequality is in fact an equality so that 
\begin{eqnarray}\label{curv-hyp}
\overline{\sigma}_0=0\quad\text{and}\quad \overline{H}'_0=\frac{R_\gamma}{2n(n-1)}.
\end{eqnarray}
On the other hand, combining (\ref{bochner}) and (\ref{ricn-1}) we get 
\begin{eqnarray*}
\overline{H}'_r-\frac{\overline{H}_r}{r}\geq 0
\end{eqnarray*}
and then the map $r\mapsto \overline{H}_r/r$ is non decreasing. This property with (\ref{curv-hyp}) gives (\ref{meancur-2}).
\qed

Recall that a compact connected Riemannian manifold is said to be of positive (respectively, negative or zero) Yamabe type when its metric can be conformally deformed into a metric with positive
(respectively, negative or zero) constant scalar curvature. This Yamabe type is precisely encoded by the sign of its Yamabe invariant ${\mathcal Y}(\partial M,[\gamma])$. In fact, this number (and so its sign) only depends on the conformal structure of the manifold so that any compact connected Riemannian manifold must belong just to one of these three conformal types. The next proposition gives a first relation between the Yamabe type of a connected component of the conformal boundary $(\partial M,[\gamma])$ and the asymptotic behavior of the isoperimetric profile of certain compact domains in $M$. More precisely, we have
\begin{proposition}\label{large-regions}
Let $(M,g)$ be an $(n+1)$-dimensional ALH manifold of order $C^{m,\alpha}$ with $m\geq 3$ and $0<\alpha<1$. Let $r$ be the geodesic defining function associated with a metric $\gamma$ in the conformal infinity $(\partial M,[\gamma])$. For each $r>0$ small enough, there exists a compact domain $\Omega_r\subset M$ with smooth boundary $\partial\Omega_r$ such that   
\begin{equation}\label{H-zero}
\lim_{r\rightarrow 0}\frac{A_r}{V_r}=n,
\end{equation}
where $A_r$ and $V_r$ denote respectively the $n$-dimensional Riemannian area of $\partial\Omega_r$ and the Riemannian volume of the domain $\Omega_r$. Moreover if we assume that a connected component of $\partial M$ has negative Yamabe invariant and that the scalar curvature $R_g$ of $(M,g)$ satisfies 
$$
R_g+n(n+1)=o(r^2)
$$
near this component for some defining function $r$, then $A_r<n\,V_r$ for all $r>0$ small enough. 
\end{proposition} 

\pf
Let $r$ be the geodesic defining function associated with a metric $\gamma$ in the conformal infinity and denote by $\Sigma_j$, $j\in\{1,\cdots,k\}$, its connected components. Let $U$ be a neighborhood of $\partial M$ and for $t>0$ sufficiently small define $M_t=M\setminus\big(r^{-1}(]0,t[)\cap U\big)$. Then there exists $t_0>0$ such that for all $0<t<t_0$, $(M_t,g)$ is a compact Riemannian manifold contained in $M$ whose boundary 
\begin{eqnarray*}
\partial M_{t}=\Sigma_{1,t}\sqcup\cdots\sqcup\Sigma_{k,t}
\end{eqnarray*}
has exactly $k$ connected components $\Sigma_{j,t}$ each of them being diffeomorphic to $\Sigma_j$ for $j\in\{1,\cdots,k\}$. Finally we fix $0<r_0<t_0$ and we define for $0<r<r_0$ the smooth compact domain $\Omega_r$ by fixing all but one the components of $\partial M_t$ at a $\overline{g}$-distance $r_0$ to $\partial M$ and, in particular, we can assume that the boundary of $\Omega_r$ is
\begin{eqnarray*}
\partial\Omega_r=\Sigma_{1,r}\sqcup\Sigma_{2,r_0}\sqcup\cdots\sqcup\Sigma_{k,r_0}.
\end{eqnarray*} 

Now by a direct application of the Taylor formula we have from (\ref{Taylor}) and (\ref{MeanCurvatureF2}) that 
\begin{eqnarray}\label{VolumeElement}
\sqrt{\frac{{\rm det\,} g_r}{{\rm det\,} \gamma}} = 1+\alpha_1r+\alpha_2 r^2+O(r^{2+\alpha})
\end{eqnarray}
where 
\begin{eqnarray}\label{alphas}
\alpha_1=-nH'_0\quad\text{and}\quad\alpha_2=\frac{n}{2}\left(nH'^2_0-\frac{H''_0}{2}\right).
\end{eqnarray} 
Recall that $H_r$ is the mean curvature of $\Sigma_{1,r}$ in $(M,g)$ for $0<r<r_0$ whose $C^{2,\alpha}$ extension to $r=0$ has been proved in the previous proposition. Then the area $A_r$ of the compact hypersurface $\partial\Omega_r$ with respect to the metric $g$ is
\begin{eqnarray*}
A_r=A(\Sigma_{1,r})+C_0= r^{-n}\int_{\Sigma_1}\sqrt{\frac{{\rm det\,} g_r}{{\rm det\,} \gamma}}\,dv_{\gamma}+C_0
\end{eqnarray*}
where $dv_\gamma$ is the Riemannian volume element of $\gamma$ and $C_0$ is the area of the other connected components of $\partial\Omega_r$ (which does not depend of $r$). A straightforward computation using (\ref{VolumeElement}) gives 
\begin{eqnarray}\label{arean3}
A_r=r^{-n}{\rm Vol}_\gamma(\Sigma_1)\left(1+\beta_1 r+\beta_2 r^2+O(r^{2+\alpha})\right)
\end{eqnarray}
for $n\geq 3$ and
\begin{eqnarray}\label{arean2}
A_r=r^{-2}{\rm Vol}_\gamma(\Sigma_1)\left(1+\beta_1 r+O(r^2)\right)
\end{eqnarray}
for $n=2$ where $\beta_j$ is the constant defined by
\begin{eqnarray}\label{betas}
\beta_j=\frac{1}{{\rm Vol}_\gamma(\Sigma_1)}\int_{\Sigma_1}\alpha_j\,dv_\gamma
\end{eqnarray}
for $j=1,2$. Here ${\rm Vol}_\gamma(\Sigma_1)$ denotes the Riemannian volume of $\Sigma_1$ with respect to $\gamma$. Similarly there exists a constant $C_1>0$ independent of $r$ such that the volume $V_r$ of $\Omega_r$ with respect to $g$ is given by
\begin{eqnarray*}
V_r= C_1+\int_{\Sigma_1}\int_r^{r_0}\sqrt{\frac{{{\rm det\,} g_s}}{{\rm det\,} \gamma}}\,dsdv_{\gamma}
\end{eqnarray*}
so that 
\begin{eqnarray}\label{volumen3}
V_r= \frac{r^{-n}{\rm Vol}_\gamma(\Sigma_1)}{n}\left(1+\frac{n\beta_1}{n-1}r+\frac{n\beta_2}{n-2}r^2+O(r^{2+\alpha})\right)
\end{eqnarray}
for $n\geq 3$ and 
\begin{eqnarray}\label{volumen2}
V_r=\frac{r^{-2}{\rm Vol}_\gamma(\Sigma_1)}{2}\left(1+2\beta_1r+2\beta_2r^2\log\frac{1}{r}+O(r^2)\right)
\end{eqnarray}
for $n=2$. Combining (\ref{arean3}) with (\ref{volumen3}) and (\ref{arean2}) with (\ref{volumen2}) immediately prove that (\ref{H-zero}) holds for all $n\geq 2$. Now if $\Sigma_1$ has negative Yamabe invariant we can assume without loss in generality that the scalar curvature of $\gamma$ is negative on $\Sigma_1$. Moreover if $R_g+n(n+1)=o(r^2)$, we have from Proposition \ref{H-limit} that $H'_0=0$ and thus 
\begin{eqnarray*}
\beta_1=0\quad\text{and}\quad\beta_2=-\frac{n}{4{\rm Vol}_\gamma(\Sigma_1)}\int_{\Sigma_1}H''_0\,dv_\gamma
\end{eqnarray*}
because of (\ref{alphas}) and (\ref{betas}). Using these facts in (\ref{arean3}) and (\ref{volumen3}) finally leads for $n\geq 3$ to
\begin{eqnarray*}
\frac{A_r}{V_r}=
n\left(1+\frac{n}{2(n-2){\rm Vol}_\gamma(\Sigma_1)}
\left(\int_{\Sigma_1}H''_0\,dv_\gamma\right)r^2+O(r^{2+\alpha})\right).
\end{eqnarray*}
Now since we assume that $R_\gamma$ is negative on $\Sigma_1$, the inequality in (\ref{meancur-1}) of Proposition \ref{H-limit} implies that 
\begin{eqnarray}\label{MeanIntegral}
\int_{\Sigma_1}H''_0\,dv_\gamma\leq\frac{1}{n(n-1)}\int_{\Sigma_1} R_\gamma\,dv_\gamma<0
\end{eqnarray}
so that $A_r<n V_r$ for $r$ sufficiently small. In a same way, for $n=2$ we derive from (\ref{arean2}) and (\ref{volumen2}) that
\begin{eqnarray*}
\frac{A_r}{V_r}=2\left(1+\frac{1}{{\rm Vol}_\gamma(\Sigma_1)}\left(\int_{\Sigma_1}H''_0\,dv_\gamma\right)r^2\log\frac{1}{r}+O(r^2)\right)
\end{eqnarray*}
which also gives that $A_r<2V_r$ for $r>0$ small enough because of (\ref{MeanIntegral}).  
\qed

From Proposition \ref{large-regions}, we immediately deduce  an upper bound for the Cheeger constant ${\mathcal Ch}(M,g)$ of any ALH manifold.  Recall that this isoperimetric constant is defined for any Riemannian manifold by 
\begin{eqnarray}\label{def-Ch}
{\mathcal Ch}(M,g)=\inf_{\Omega}\frac{A(\partial\Omega)}{V(\Omega)},
\end{eqnarray}
where the infimum is taken over all the compact (smooth) domains in $M$, $A(\partial\Omega)$ being the area of the compact hypersurface $\partial \Omega$ and $V(\Omega)$ the volume of the domain $\Omega$. Then we prove
\begin{corollary}\label{upper-Cheeger}
The Cheeger isoperimetric constant ${\mathcal Ch}(M,g)$ of an $(n+1)$-dimensional ALH manifold $(M,g)$ of order $C^{m,\alpha}$ with $m\geq 3$ and $0<\alpha<1$ satisfies
$$
{\mathcal Ch}(M,g)\le n.
$$
If some connected component of the conformal infinity $(\partial M,[\gamma])$ has negative Yamabe invariant and if, near  this component, for some defining function $r$, the scalar curvature $R_g$ of $M$ satisfies
$$
R_g+n(n+1)=o(r^2),
$$
 then 
$$
{\mathcal Ch}(M,g)<n.
$$
\end{corollary}

\begin{remark}
{\rm
Note that the second part of the previous corollary applies for WPE manifolds.}
\end{remark}

\pf 
Choose any geodesic defining function $r$ for the ALH manifold $M$ and let $\Omega_r\subset M$ be the compact domain associated with $r>0$ small enough as in Proposition \ref{large-regions}. If ${\mathcal Ch}(M,g) > n$, since Proposition \ref{large-regions} assures that
$\lim_{r\rightarrow 0}(A(\partial\Omega_r)/V(\Omega_r))=n$, then ${\mathcal Ch}(M,g)$ could not be a lower bound for the set of the isoperimetric quotients $A(\partial\Omega)/V(\Omega)$, with $\Omega
\subset M$ compact. Thus, we have ${\mathcal Ch}(M,g)\le n$ for any ALH manifold $M$. Now, assume that the conformal infinity $(\partial M,[\gamma])$ of the given ALH manifold $M$ has at least one connected component with negative Yamabe type which, near this component, satisfies the decay condition $R_g+n(n+1)=o(r^2)$. In this situation, the second part of Proposition \ref{large-regions} provides compact domains $\Omega_r$ in $M$ such that $A(\partial\Omega_r)<n\,V(\Omega_r)$, hence we directly have ${\mathcal Ch}(M,g)<n$.
\qed


\section{Some examples}


Let $\big(B^{n+1},\frac{4|dx|^2}{(1-|x|^2)^2}\big)$ be the Poincar{\'e} hyperbolic ball. Using the change of variables given by $s=\ln\frac{1+|x|}{1-|x|}\in{\mathbb R}_+$, we obtain
the metric 
$$
g=ds^2+(\sinh^2 s)\gamma_{{\mathbb S}^n}.
$$
This expression for the Poincar{\'e} metric is valid only on the punctured ball $B^{n+1}-\{0\}\cong {\mathbb R}_+^*\times {\mathbb S}^n$, although it is smoothly extendible to the origin. 

Written in this way, we can see that the hyperbolic metric is an example of the so-called {\em warped Riemannian products} (see for instance, \cite{Be,O'N,K}). In general, if $I\subset {\mathbb R}$ is an open interval, $(P,\gamma)$ a Riemannian $n$-manifold and $f\in C^\infty(I)$ is a positive function, we will say that the $(n+1)$-dimensional Riemannian manifold $(I\times P, g=ds^2+f(s)^2 \gamma)$ is the  product of $I$ and $P$  warped by 
the function $f$.  We will restrict ourselves to warping functions
$f$ satisfying $f''-f=0$. With this choice, we ensure  that $\Ric_g(\frac{\partial}{\partial s},\frac{\partial}{\partial s})=-n$ at each point of $I\times P$ (see \cite[Lemma 4]{K}). Moreover, taking also into account the values of $\Ric_g$ on the directions orthogonal to the vector field $\frac{\partial}{\partial s}$, that is, directions tangent to $P$, we conclude that there are essentially three types of warped products which eventually may produce PE spaces with scalar curvature $-n(n+1)$, according to the
warping function is chosen to be $\sinh s$, $e^s$ or $\cosh s$ (conical singularities and cusps being provisionally permitted). 

\begin{example}\label{rk3}
{\rm
The first class of manifolds we consider here is the so-called {\em hyperbolic cones} on given compact Riemannian manifolds $(P,\gamma)$, which are defined by
$$
\big(M={\mathbb R}_+\times P, g = ds^2 + (\sinh^2 s) \gamma\big).
$$  
Defining a new variable $t\in ]0,1]$ by  $t=\tanh\frac{s}{2}$, we obtain that the conformal metric
$$
{\overline g}=\left(\frac1{1+\cosh s}\right)^2 g = dt^2 + t^2 \gamma
$$
extends to $[0,1]\times P$ that is to a compact manifold with boundary $\{1\}\times P\cong P$ and a
conical singularity at $t=0$. This singularity is removable if and only if $(P, \gamma)$ is the round unit $n$-sphere and, in this case, the corresponding hyperbolic cone is nothing but the $(n+1)$-dimensional hyperbolic space (see \cite[p. 269, Lemma 9.114]{Be}). 

It is clear from the above considerations,  that if $({\mathbb S}^n,\gamma)$ is the round unit sphere and $f:{\mathbb R}_+\rightarrow{\mathbb R}^*_+$ coincides with $s\mapsto\sinh s$, both near $0$ and $+\infty$, the warped product metric given by 
$$
\big(M={\mathbb R}_+\times {\mathbb S}^n, g =ds^2 + f(s)^2 \gamma\big),
$$
  also avoids the conical singularity at $s=0$ and that the resultant $(n+1)$-dimensional Riemannian manifold $(M,g)$ is a rotationally invariant deformation of the hyperbolic space. These manifolds are not in general PE spaces but are WPE spaces whose conformal infinity is obviously $({\mathbb S}^n, [\gamma])$, and so they trivially have positive Yamabe type. For a fixed $\varepsilon>0$, if we choose $f$ in such a way that
$$
f \left(\frac1{\varepsilon}\right) = f \left(\frac{3}{\varepsilon}\right)=\varepsilon^{\frac1{n}},
\qquad f(s)\ge \varepsilon^{\frac1{n}},\quad\forall s\in\left[\frac1{\varepsilon},\frac{3}{\varepsilon}\right],
$$    
 and represent by $\Omega_\varepsilon$ the domain  $\left]\frac1{\varepsilon},\frac{3}{\varepsilon}\right[
\times P\subset M$, we have
$$
\frac{A(\partial\Omega_\varepsilon)}{V(\Omega_\varepsilon)}=\frac{2\varepsilon A(P)}{\int_{1/\varepsilon}
^{3/\varepsilon}\int_Pf(s)^n\,ds\,dP}\le \varepsilon.
$$
Hence, from Definition (\ref{def-Ch}), we obtain
$$
{\mathcal Ch}(M,g)\le \varepsilon.
$$
This means that we have  $(n+1)$-dimensional ALH manifolds (in fact, WPE spaces) with conformal infinity of positive Yamabe type and Cheeger isoperimetric constants arbitrarily small belonging to the interval $[0,n]$.}
\end{example}

\begin{example}{\rm
The second type we also consider here is of the form
$$
\big(M={\mathbb R}\times P,  g = ds^2 + (\cosh^2 s) \gamma\big).
$$
Such warped product metrics satisfy (see \cite{K}, for instance)
 $$\Ric_g(\frac{\partial}{\partial s},\frac{\partial}{\partial s})=-n \quad {\rm and}  \quad
R_g+n(n+1)=R_\gamma+n(n-1).
$$
Hence, as before, to get a WPE space, we will impose on $(P,\gamma)$ to have scalar curvature $-n(n-1)$.  In order to compactify, we define a new variable $t\in]0,\pi[$ by the relation $t=2\arctan e^s$. Then
$$
{\overline g}=\left(\frac1{\cosh s}\right)^2g=dt^2+\gamma,
$$
is smoothly extendable to the compact manifold $[0,\pi]\times P$. Hence, its conformal infinity $(\partial M, [\gamma])$ consists of two copies of $(P,[\gamma])$. Thus we obtain an example of WPE space whose conformal compactification has two connected components at the conformal boundary (a {\em wormhole} in the physical jargon), both with negative Yamabe type. 
Now, if $f:{\mathbb R}\rightarrow{\mathbb R}^*_+$ coincides with $s\mapsto\cosh s$ both near $-\infty$ and $+\infty$, the manifold $(P, \gamma)$ is a Riemannian manifold with constant scalar curvature $-n(n-1)$, and we consider a warped product 
$$
\big(M={\mathbb R}\times P, g = ds^2 + f(s)^2 \gamma\big),
$$
 then  the corresponding $(n+1)$-dimensional Riemannian manifold $(M, g)$ is a deformation of the above one which is still a WPE space with conformal infinity consisting of two copies
of $(P, [\gamma])$ as well. We may require $f$ to have exactly the same behavior as in the above example on the interval $[\frac1{\varepsilon},\frac{3}{\varepsilon}]$ and we conclude that there are $(n+1)$-dimensional ALH  manifolds (in fact, WPE spaces) with conformal infinity of negative Yamabe type and Cheeger isoperimetric constants arbitrarily small inside the interval $[0,n[$.}
\end{example}


\section{The Cheeger constant of $(M, g)$ and the Yamabe type of $(\partial M, [\gamma])$}


In this section, we state and prove the main result of this paper which relates the value of the Cheeger constant of a conformally compact Riemannian manifold with the Yamabe type of its conformal infinity. As we are admitting the possibility that $\partial M$ is not connected, we should be more precise in the definition of the Yamabe type of the conformal infinity in this setting. However, using a famous result on the connectedness of the boundary at infinity by E. Witten and S.-T. Yau (for boundaries with a component of positive Yamabe type) and by M. Cai and G. Galloway (for boundaries with a component of null Yamabe type), we will avoid this discussion. In order to make this paper self-contained, we include a new proof of this result which simplifies and unifies the two aforementioned proofs and slightly weakens their hypotheses (and fits into ours). In fact, we will show that these two theorems can be seen as direct consequences of an old paper  by A. Kasue \cite{Ka} generalizing the Bonnet-Myers theorem to complete manifolds with non-empty boundary. 
\begin{theorem}{\rm (\cite{WY,CG})}\label{connectedness} 
Let $(M,g)$ be an $(n+1)$-dimensional conformally compact Riemannian manifold of order $C^{m,\alpha}$ with $m\geq 3$ and $0<\alpha<1$ whose Ricci tensor and scalar curvature satisfy
\begin{equation}\label{Ricci-Scalar}
\hbox{\rm Ric\,}_g+ng\ge 0,\qquad R_g+n(n+1)=o(r^2),
\end{equation}
where $r$ is a geodesic defining function. Suppose that the conformal infinity $(\partial M,[\gamma])$ has a connected component with non-negative Yamabe invariant. Then
$\partial M$ is connected.
\end{theorem}

\pf 
Let $\partial M_0$ be the component of the conformal infinity $(\partial M,[\gamma])$ with non-negative Yamabe invariant. By the solution of the Yamabe problem \cite{Sc}, we can choose a  metric $\gamma\in [\gamma]$ with constant scalar curvature, say $R_{\gamma}=n(n-1)\varepsilon^2$, where $\varepsilon=1$ or $\varepsilon=0$ depending on whether the Yamabe invariant is positive or zero. Denote by $r$ the unique geodesic defining function associated with $\gamma$ whose existence is given by Lemma \ref{defining-geodesic}. Let $U$ be a neighborhood of $\partial M_0$ which does not meet any other component of $\partial M$ and $M_t=M-\big(r^{-1}(]0,t[)\cap U\big)$. For $t>0$ sufficiently small,
$(M_t,g)$ is a Riemannian manifold with boundary, whose one of its components $\Sigma_t$ is diffeomorphic to $\partial M_0$ and whose mean curvature satisfies  
\begin{eqnarray}\label{mean-cur3}
H_t\ge  1+\frac{\varepsilon^2}{2}t^2,
\end{eqnarray}
(see (\ref{meancur-2}) in Proposition \ref{H-limit}).

Suppose now that the Yamabe invariant of $\partial M_0$ is positive, that is, $\varepsilon=1$. In this situation, $(M_t,g)$ is a complete Riemannian manifold such that $\hbox{\rm Ric\,}_g\geq -ng$ and whose compact connected boundary $\Sigma_t$ satisfies $H_t>1$ from (\ref{mean-cur3}). Hence we can apply \cite[Theorem A]{Ka} (see also \cite[Proposition 2]{CG}, which reproves a part of Kasue's result) and conclude that $M_t$ must be compact. So $\Sigma_t$ is the unique component of its boundary. This means that $\partial M_0$ is the unique component in $\partial M$ and $\partial M$ is connected, as claimed.

When the Yamabe invariant of $\partial M_0$ is zero, i.e. $\varepsilon=0$, the same reasoning provides $(M_t,g)$ as above with $H_t\ge 1$. If $M_t$ is non-compact, a direct application of \cite[Theorem C]{Ka} implies that $\min H_t= 1$ and $M_t$ is isometric to the warped product $\big([0,+\infty[\times \Sigma_t,ds^2+e^{-2s}g_{|\Sigma_t}\big)$ (for definitions and properties about warped products, see \cite{Be,O'N,K} or Remark \ref{rk3}  below). This means that $M_t$ has a connected compact boundary $\Sigma_t$ and one end which is a {\em hyperbolic cusp}. This contradicts the fact that every end of an ALH manifold has infinite volume (since the volume of a hyperbolic cusp is finite). We conclude that $M_t$ is compact and then $\partial M$ has to be connected.
\qed

Assume now that $(M,g)$ is a conformally compact manifold satisfying the curvature assumptions (\ref{Ricci-Scalar}) and that a connected component of its conformal boundary $(\partial M,[\gamma])$ has negative Yamabe invariant. From the solution of the Yamabe problem, we can assume that this component has constant negative scalar curvature so that Corollary \ref{upper-Cheeger} implies that ${\mathcal Ch}(M,g)<n$. Since ${\mathcal Ch}(M,g)\leq n$ for all ALH manifolds, we immediately deduce that if ${\mathcal Ch}(M,g)= n$, then a connected component of the boundary at infinity $(\partial M,[\gamma])$ must have a non negative Yamabe invariant and this implies that $\partial M$ is in fact connected using Theorem \ref{connectedness}. We sum up these properties and show that the converse is also true in our main result:
\begin{theorem}\label{maintheorem}
The Cheeger constant ${\mathcal Ch}(M,g)$ of an $(n+1)$-dimensional ALH manifold $(M,g)$ of order $C^{m,\alpha}$ with $m\geq 3$ and $0<\alpha<1$ whose conformal infinity is $(\partial M,[\gamma])$ satisfies
$$
{\mathcal Ch}(M,g)\le n.
$$
Moreover, if the Ricci and the scalar curvatures of $M$ satisfy 
$$
\hbox{\rm Ric}_g+ng\ge 0,\qquad R_g+n(n+1)=o\big(r^2\big),
$$
where $r$ is any defining function on $M$, then we have
$$
{\mathcal Ch}(M,g)= n\quad\Longleftrightarrow\quad{\mathcal Y}(\partial M,[\gamma])\ge 0,
$$
where ${\mathcal Y}(\partial M,[\gamma])$ denotes the Yamabe invariant of the conformal boundary. In particular, $(\partial M,[\gamma])$ is connected.
\end{theorem}

\pf
It only remains to prove that if ${\mathcal Y}(\partial M,[\gamma])\ge 0$ then ${\mathcal Ch}(M,g)= n$. Since we know from Corollary \ref{upper-Cheeger} that ${\mathcal Ch}(M,g)\leq n$, hence it is sufficient to prove that ${\mathcal Ch}(M,g)\geq n$.   We originally proved this fact using an approach relying on geometric measure theory (as in \cite{W}). However, we found in \cite{FR} an elementary proof due to Gilles Carron for PE manifolds. We shall see that these arguments always work under our (weaker) assumptions. Indeed, recall from \cite{L} (see \cite[Proposition 3]{HM} for a precise statement) that if $r$ is a geodesic defining function satisfying  (\ref{Ricci-Scalar}), then there exists a unique positive function $u\in C^\infty (M)$ such that $\Delta u=(n+1)u$ and $u-\frac1{r}$ is a bounded function. If in addition ${\mathcal Y}(\partial M, [\gamma])\geq 0$, it can be shown (see \cite[Proposition 4.2]{L}) that the function $u^2-|\nabla u|_g^2$ is superharmonic on $M$, extends continuously to $\overline{M}$ and vanishes on $\partial M$ so that the strong minimum principle implies $u^2-|\nabla u|_g^2\geq 0$. Then if we let $f=\ln u$ on $M$, it is straightforward to check that this function satisfies 
\begin{eqnarray}\label{SuperMeanExit}
|\nabla f|_g^2\leq 1\quad\text{and}\quad\Delta f\geq n. 
\end{eqnarray}
Now consider a smooth compact domain $\Omega$ in $M$ and integrate (using the Stokes formula) the second inequality in (\ref{SuperMeanExit}) over $\Omega$, to get 
\begin{eqnarray}\label{In1}
-\int_{\partial\Omega}\frac{\partial f}{\partial N}\geq n V(\Omega)
\end{eqnarray}
where $N$ denotes the inner unit normal to $\Sigma$ in $\Omega$. On the other hand, the first inequality in (\ref{SuperMeanExit}) implies 
\begin{eqnarray}\label{In2}
-\int_{\partial\Omega}\frac{\partial f}{\partial N}\leq
\int_{\partial\Omega}\Big|\frac{\partial f}{\partial N}\Big|\leq 
A(\partial\Omega).
\end{eqnarray}
Combining (\ref{In1}) and (\ref{In2}) gives  $A(\partial\Omega)\geq n V(\Omega)$ for all smooth compact domains of $M$ which is exactly ${\mathcal Ch}(M,g)\geq n$, as claimed.  
\qed 

\begin{remark}\label{r1}{\rm
Theorem \ref{maintheorem}  implies that the rotationally invariant deformations $(M,g)$ of the hyperbolic space
${\mathbb H}^{n+1}$ built in Example \ref{rk3} to provide examples of WPE spaces with conformal infinity of positive Yamabe type and arbitrarily small Cheeger constant, cannot satisfy the hypothesis on the Ricci curvature, that is, $\hbox{\rm Ric}_g$ cannot admit $-n$ as a lower bound. If so, these examples provide conformally compact manifolds whose conformal infinity is the round conformal sphere and such that (\ref{Ricci-Scalar}) holds. However, the rigidity result of such conformally compact manifolds implies that $(M,g)$ has to be the hyperbolic space that is $f(s)=\sinh s$ for all $s\in\mathbb{R}_+$ (see Corollary $1.5$ in \cite{LQS}). This contradicts the fact that the Cheeger constant of the hyperbolic space is $n$.  It is worth mentioning that this result can be directly (and easily) observed by looking closer at such examples. Indeed, if $(M={\mathbb R}_+\times{\mathbb S}^n,g=ds^2+f(s)^2\gamma)$ is as in Example \ref{rk3} with $\hbox{\rm Ric}_g+ng\ge 0$, we would necessarily have $f''\le f$. Letting $y=f'/f$, we immediately observe that $y$ satisfies $y'+y^2\leq 1$ on $\mathbb{R}_+$. Since $y$ coincides with $z : s\mapsto{\rm cotanh\,}s$ in a neighborhood of $0$ and $+\infty$ satisfying $z'+z^2=1$, we can apply Lemma $4.1$ in \cite{ballmann} to conclude first that $y(s)={\rm cotanh\,} s$ for all $s\in\mathbb{R}_+$ and then that $f(s)=\sinh s$ on $\mathbb{R}_+$: the manifold $(M,g)$ is isometric to the hyperbolic space so that we get the desired contradiction.   
}
\end{remark}


\section{Some direct consequences}


\subsection{Minimizer of the Cheeger constant} A direct consequence of Theorem \ref{maintheorem} is that the Cheeger constant of a conformally compact manifold satisfying (\ref{Ricci-Scalar}) does not possess smooth minimizer. Indeed if we denote by $\Omega_0$ such a domain then we obviously have 
\begin{eqnarray}\label{minimizers}
{\mathcal Ch}(M,g)=n=\frac{A(\partial\Omega_0)}{V(\Omega_0)}.
\end{eqnarray}
On the other hand, given a smooth function $f$ on $\partial\Omega_0$, we consider the normal variation of $\partial\Omega_0$ defined by 
\begin{eqnarray*}
\psi_t:p\in\partial\Omega_0\mapsto \exp_p\big(-tf(p)N_0(p)\big)\in M
\end{eqnarray*}
where $\exp$ is the exponential map of $M$ and $N_0$ is the inner unit vector normal to $\partial\Omega_0$ in $\Omega_0$. Denote by $A(t)$ the area of the hypersurface $\psi_t(\partial\Omega_0)$ as well as $V(t)$ the volume of the domain enclosed by $\psi_t(\partial\Omega_0)$ for $|t|$ sufficiently small. Since $\Omega_0$ is a minimum of ${\mathcal Ch}(M,g)$ we must have
\begin{eqnarray*}
\frac{d}{dt}_{|t=0}\frac{A(t)}{V(t)}=0
\end{eqnarray*}
which, from the first variational formulae of the area and of the volume, is equivalent to 
\begin{eqnarray*}
\int_{\partial\Omega_0}f\Big(nHV(\Omega_0)-A(\partial\Omega_0)\Big)=0
\end{eqnarray*}
for any $f\in C^\infty(\partial\Omega_0)$. We conclude that $nH=A(\partial\Omega_0)/V(\Omega_0)$ and then $H=1$ because of (\ref{minimizers}). Now since $\hbox{\rm Ric}_g+ng\ge 0$, the Heintze-Karcher inequality \cite{HK} implies 
\begin{eqnarray*}
V(\Omega_0)\leq A(\partial\Omega_0)\int_0^{R_0}\big(\cosh t-H_0\sinh t\big)^ndt
\end{eqnarray*}
where $H_0$ is the minimum of $H$ on $\partial\Omega_0$ and $R_0$ is the inradius of $\Omega_0$. As $H_0=H=1$ we immediately deduce  
\begin{eqnarray*}
nV(\Omega_0)\leq(1-e^{-nR_0})A(\partial\Omega_0)<A(\partial\Omega_0)
\end{eqnarray*} 
and this precisely contradicts (\ref{minimizers}).

\subsection{Isoperimetric inequality} Note that, by the very definition of the Cheeger constant, if ${\mathcal Ch}(M,g)=n$, we have that the isoperimetric inequality (\ref{LII}) is satisfied on $M$. Conversely, when this inequality holds for each compact domain in $M$, we can only  conclude the inequality ${\mathcal Ch}(M,g)\ge n$. But, if $M$ is an ALH manifold, from the first inequality in Corollary \ref{upper-Cheeger}, we have that the corresponding equality is achieved. Thus, we obtain another characterization of the non-negativity of the conformal infinity of the class of ALH manifolds that we are studying.  
\begin{corollary}\label{yamabe-yau}
Let $M$ be an $(n+1)$-dimensional conformally compact manifold of order $C^{m,\alpha}$ with $m\geq 3$ and $0<\alpha<1$ whose conformal infinity is $(\partial M,[\gamma])$.  If its Ricci and scalar curvatures satisfy (\ref{Ricci-Scalar}) then we have ${\mathcal Y}(\partial M,[\gamma])\ge 0$ if and only if the isoperimetric inequality
$$
A(\partial\Omega)\ge n V(\Omega),
$$
holds for any compact domain $\Omega\subset M$.
\end{corollary} 

In the particular case where the ALH manifold is a hyperbolic manifold, taking into account the fundamental work \cite[Theorem 4.7] {SY} by R. Schoen and S.-T. Yau about conformally flat manifolds, Theorem \ref{maintheorem} allows to decide when the linear hyperbolic isoperimetric inequality remains  valid after quotienting by a Kleinian group.

\begin{corollary}
Let ${\mathbb H}^{n+1}/\Gamma$ be a geometrically finite and cusp-free quotient of the $(n+1)$-dimensional hyperbolic space by a Kleinian group. Denote by $\Lambda (\Gamma)\subset {\mathbb S}^n$ the limit set of $\Gamma$ and by ${\mathcal H}\big(\Lambda(\Gamma)
\big)$ its Hausdorff measure. Then we have$$
A(\partial\Omega)\ge n V(\Omega),\quad\forall\Omega\subset{\mathbb H}^{n+1}/\Gamma 
\quad
\Longleftrightarrow 
\quad
{\mathcal H}\big(\Lambda(\Gamma)\big)\le \frac{n-2}{2}.$$
\end{corollary}  

\subsection{Principal eigenvalue of the $p$-Laplacian} Another immediate consequence of Theorem \ref{maintheorem} is an extension of the result obtained by J. Lee  \cite{L} on the infimum of the $L^2$ spectrum of the Laplacian to the case of the $p$-Laplacian. We first briefly recall some well-known facts on this operator and its principal eigenvalue (for more details we refer to \cite{Mat,SW} and references therein).

On a Riemannian manifold, for $1<p<\infty$ and $u\in C^\infty(M)$, the $p$-Laplacian $\Delta_p$ is defined  by
\begin{eqnarray*}
\Delta_p u={\rm div\,}(|\nabla u|_g^{p-2}\nabla u).
\end{eqnarray*}
 The principal eigenvalue $\lambda_{1,p}(M)$ of the $p$-Laplacian is the maximum constant $\lambda$ such that the equation
\begin{eqnarray*}
\Delta_p u =-\lambda\, u^{p-1}
\end{eqnarray*}
admits a positive solution. Alternatively, it may be characterized variationally as the best constant in the inequality 
\begin{eqnarray*}
\lambda_{1,p}(M)\int_M |v|^p\leq\int_M|\nabla v|_g^p
\end{eqnarray*}
for any smooth and compactly supported function $v$ on $M$. From \cite{SW} we know that if $(M,g)$ is a complete $(n+1)$-dimensional Riemannian manifold $(M,g)$ with $\hbox{\rm Ric}_g+ng\ge 0$, this principal eigenvalue satisfies an inequality analogous to the famous Cheng inequality for the bottom of the $L^2$ spectrum of the standard Laplacian, namely
\begin{eqnarray}\label{Cheng}
\lambda_{1,p}(M)\leq\Big(\frac{n}{p}\Big)^p.
\end{eqnarray}
On the other hand, we claim that the Cheeger-type inequality 
\begin{eqnarray}\label{p-Cheeger}
\lambda_{1,p}(M)\geq\Big(\frac{{\mathcal Ch}(M,g)}{p}\Big)^p
\end{eqnarray}
is also satisfied. Indeed, first note that if $(\Omega_i)$ is an exhaustion of $M$ by compact domains, it is straightforward to show that
\begin{eqnarray}\label{LimSp}
\lambda_{1,p}(M)=\lim_{i\rightarrow\infty}\lambda_{1,p}^D(\Omega_i)
\end{eqnarray}
where $\lambda_{1,p}^D(\Omega_i)$ is the first eigenvalue of the $p$-Laplacian on $\Omega_i$ with Dirichlet boundary condition that is 
$$
\left\lbrace
\begin{array}{l}
\Delta_p v_i=-\lambda_{1,p}^D(\Omega_i) v_i \\
v_{i|\partial\Omega_i}=0.
\end{array}
\right.
$$
Moreover, it is proved in \cite{T} that
\begin{eqnarray}\label{p-CheegIn}
\lambda_{1,p}^D(\Omega_i)\geq\Big(\frac{{\mathcal Ch}(\Omega_i)}{p}\Big)^p
\end{eqnarray} 
where ${\mathcal Ch}(\Omega_i)$ is the Cheeger constant of $\Omega_i$ defined by
\begin{eqnarray*}
{\mathcal Ch}(\Omega_i)=\inf_{\Omega}\frac{A(\partial\Omega)}{V(\Omega)}
\end{eqnarray*}
where $\Omega$ ranges over all smooth compact domain in $\Omega_i$ and smooth boundary $\partial\Omega$. From the definition (\ref{def-Ch}) of the Cheeger constant of $M$ it is obvious that ${\mathcal Ch}(\Omega_i)\geq{\mathcal Ch}(M,g)$ for all $i$ so that (\ref{p-Cheeger}) follows directly from (\ref{LimSp}) and (\ref{p-CheegIn}). Finally applying Theorem \ref{maintheorem} to (\ref{p-Cheeger}) and combining with (\ref{Cheng}) indeed leads to the aforementioned generalization of Lee spectral estimate:
\begin{theorem}\label{lambda}
Let $(M,g)$ be an $(n+1)$-dimensional conformally compact manifold of order $C^{m,\alpha}$ with $m\geq 3$ and $0<\alpha<1$ whose conformal infinity is $(\partial M,[\gamma])$. Suppose that its Ricci and scalar curvatures  satisfy (\ref{Ricci-Scalar}). For $1<p<\infty$, if ${\mathcal Y}(\partial M,[\gamma])\ge 0$, then  
$$
\lambda_{1,p}(M)= \Big(\frac{n}{p}\Big)^p,
$$
where $\lambda_{1,p}(M)$ denotes the principal eigenvalue of the $p$-Laplacian of $M$.
\end{theorem}

\begin{remark}{\rm
According to Theorem \ref{lambda}, all the conformally compact manifolds $(M,g)$ satisfying $\hbox{Ric}_g+ng\ge 0$, a second order scalar curvature decay and ${\mathcal Y}(\partial M,[\gamma])\ge 0$, are examples of complete Riemannian manifolds with optimal Cheeger inequality (\ref{p-Cheeger}). Instead, the results obtained by D. Sullivan in \cite{Su} and by R. Schoen and S.-T. Yau in \cite{SY} imply that the geometrically finite and cusp free quotients ${\mathbb H}^{n+1}/\Gamma$ with $\frac{n-2}{2}<{\mathcal H}\big(\Lambda(\Gamma)\big)\le\frac{n}{2}$ have $\lambda_{1,2}({\mathbb H}^{n+1}/\Gamma)=\frac{n^2}{4}$ and ${\mathcal Y}({\mathbb H}^{n+1}/\Gamma)<0$. Then, we deduce from our Theorem \ref{maintheorem} that ${\mathcal Ch}({\mathbb H}^
{n+1}/\Gamma)<n$. Thus these hyperbolic manifolds give examples of PE spaces where the Cheeger
inequality is not sharp.            
}
\end{remark}


\begin{thebibliography}{BHHMm}
\bibitem [Ba]{ballmann} W. Ballmann, {\em Riccati equation and volume estimates}, preprint 2016.

\bibitem [Be]{Be} A. Besse, {\em Einstein manifolds}, Springer, New York, (1987).

\bibitem [Bi]{Bi} O. Biquard (ed.), {\em AdS/CFT 
Correspondence: Einstein Metrics and Their Conformal Boundaries}, 
IRMA Lectures in Mathematics and Theoretical
Physics, Euro. Math. Soc., Z{\"u}rich, 2005.

\bibitem [BZ]{BZ} Y. D. Burago, V. A. Zalgaller, {\em Geometric inequalities}, Springer-Verlag, Berlin / Heidelberg, 1988. 

\bibitem [CG]{CG} M. Cai, G. Galloway, {\em Boundaries of zero scalar curvature in the 
AdS/CFT correspondence}, Adv. Theor. Math. Phys. {\bf 3} (1999), 1769--1783.

\bibitem [C]{C} J. Cao, {\em Cheeger isoperimetric constants of Gromov-hyperbolic 
spaces with quasi-poles}, Commun, Contemp. Math. {\bf 4} (2000), 511--533.

\bibitem [Ch]{Ch} J. Cheeger, {\em A lower bound for the smallest eigenvalue of the Laplacian}, Problems in Analysis, 195--199, Princeton University Press, 1970.

\bibitem [CDLS]{CDLS} P. T. Chru{\'s}ciel, E. Delay, J. M. Lee, D. N. Skinner, {\em Boundary regularity of conformally compact Einstein metrics}, J. Diff. Geom. {\bf 69} (2005), 111--136.

\bibitem [E]{E} J. F. Escobar, {\em Conformal deformation of a Riemannian metric to a scalar 
flat metric with constant mean curvature on the boundary}, Ann. of Math., {\bf 136} (1992), 1--50.  


\bibitem [FG]{FG} C. Fefferman, C. R. Graham, {\em Conformal invariants}, in ``{\'E}lie Cartan
et les math{\' e}matiques d'aujourd'hui", Ast{\' e}risque {\bf } (1985), 95--116.

\bibitem [FR]{FR} F. Ferrari, A. Rovai, {\em Holography, probe branes and isoperimetric inequalities}, Phys. Letters B {\bf 747} (2015), 212--216. 
\bibitem [GL]{GL} C. R. Graham, J. M. Lee, {\em Einstein metrics with prescribed conformal 
infinity on the ball}, Adv. Math. {\bf 87} (1991), 186--225. 

\bibitem [GQ]{GQ} C. Guillarmou, J. Qing, {\em Spectral characterization of Poincar{\'e}-Einstein
manifolds with infinity of positive Yamabe type}, Int. Math. Res. Not. {\bf 28} (2010), 1720-1740.

\bibitem [HK]{HK} E. Heintze, H. Karcher, {\em A general comparison theorem with applications to volume estimates for submanifolds}, Ann. Sci. \'Ecole Norm. Sup. {\bf 11} no. 4 (1978), 451--470.

\bibitem [HM]{HM} O. Hijazi, S. Montiel, {\em Supersymmetric rigidity of asymptotically locally hyperbolic manifolds}, Int. J. Math. {\bf 25} no. 3 (2014), 25 pages.

\bibitem [K]{K} M. Kanai, {\em On a differential equation characterizing a Riemannian structure of a manifold}, Tokyo J. Math. {\bf 6} (1983), 143--151.

\bibitem [Ka]{Ka} A. Kasue, {\em Ricci curvature, geodesics and some geometric properties of Riemannian manifolds with boundary}, J. Math. Soc. Japan {\bf 35} (1983), 117--131.

\bibitem [L]{L} J. M. Lee, {\em The spectrum of an asymptotically hyperbolic Einstein manifold}, Comm. 
Anal. Geom. {\bf 2} (1995), 253--271.

\bibitem [LQS]{LQS} G. Li, J. Qing, Y. Shi {\em Gap phenomena and curvature estimates for conformally compact Einstein manifolds}, Trans. Amer. Math. Soc. {\bf 369}  no.6 (2017), 4385--4413.

\bibitem [Mat]{Mat} A.-M. Matei, {\it First eigenvalue for the $p$-Laplace operator}, Nonlinear Analysis {\bf 39} (2000), 1051-1068.

\bibitem [Maz]{M} R. Mazzeo {\it Unique continuation at infinity and embedded eigenvalues for asymptotically hyperbolic manifolds}, Amer. J. Math. {\bf 113} (1991), 25-45.

\bibitem [MP]{MP} R. Mazzeo, F. Pacard, {\it Constant curvature foliations in asymptotically hyperbolic spaces}, Rev. Mat. Iberoam. {\bf 27} (2011), 303-333.
 
\bibitem [O'N]{O'N} B. O'Neill, {\em Semi-Riemannian Geometry}, Academic Press,
1983.

\bibitem [P]{P} P. Petersen, {\em Riemannian Geometry}, Graduate Texts in Mathematics, vol. 171, 2nd edn. Springer, New York (2006).

\bibitem [Q]{Q} J. Qing, {\em On the rigidity for conformally compact Einstein manifolds}, Intern. 
Math. Res. Not. {\bf 21} (2003), 1141--1153.

\bibitem [Sc2]{Sc} R. Schoen, {\em Conformal deformation of a Riemannian metric to constant
scalar curvature}, J. Diff. Geom. {\bf 20} (1984), 479--495.

\bibitem [SY]{SY} R. Schoen, S.-T. Yau, {\em Conformally flat manifolds, Kleinian groups and scalar curvature}, Invent. Math. {\bf 92} (1988), 47--71.

\bibitem [Su]{Su} D. Sullivan, {\em Related aspects of positivity in Riemannian geometry}, J. Diff. Geom. {\bf 25} (1987), 327--351.

\bibitem [SW]{SW} C.-J.A. Sung, J. Wang, {\em Sharp gradient estimate and spectral rigidity for $p$-Laplacian}, Math. Res. Lett. {\bf 21} (4) (2014), 885--904.

\bibitem [T]{T} H. Takeuchi, {On the first eigenvalue of the $p$-Laplacian in a Riemannian manifold}, Tokyo J. Math. {\bf 21} (1) (1998), 135--140.

\bibitem [W]{W} X. Wang, {\em A new proof of Lee's theorem on the spectrum of conformally compact Einstein manifold}, Comm. Anal. Geom. {\bf 10} (2002), 647--651.

\bibitem [WY]{WY} E. Witten, S.-T. Yau, {\em Connectedness of the boundary in the AdS/CFT correspondence}, Adv. Theor. Math. Phys. {\bf 3} (1999), 1635--1655. 

\bibitem [Y]{Y} S.-T. Yau, {\em Isoperimetric constants and the first eigenvalue of a compact Riemannian manifold}, Ann. Sci. {\'E}cole Norm. Sup. {\bf 8} (1975), 487--507.

\end{thebibliography}
\end{document}